\documentclass[conference]{IEEEtran}
\IEEEoverridecommandlockouts
\usepackage{cite}
\usepackage{amsmath,amssymb,amsfonts}
\usepackage{textcomp}

\usepackage{graphicx} 
\usepackage{booktabs} 
\usepackage{subcaption}
\usepackage{subfiles}
\usepackage{url}
\usepackage{amssymb}
\usepackage{amsmath}
\usepackage{nicefrac}
\usepackage{multirow}

\usepackage{algorithm}
\usepackage[noend]{algpseudocode}
\usepackage{hyperref}

\usepackage{xcolor}

\usepackage{pifont}

\usepackage{amsmath,mathrsfs,amssymb,mathtools,bm}
\usepackage{scalerel}
\usepackage{amsthm}
\usepackage{txfonts}
\usepackage{fontawesome}
\usepackage{wrapfig}
\usepackage{todonotes}
\usepackage{xcolor}
\usepackage{setspace}
\usepackage{bm,bbm}
\usepackage{paralist}
\usepackage{enumitem}

 \newcommand{\cF}{\mathcal{F}}

\usepackage{esvect}

\makeatletter
\newcommand{\oset}[3][0ex]{%
  \mathrel{\mathop{#3}\limits^{
    \vbox to#1{\kern-2\ex@
    \hbox{$\scriptstyle#2$}\vss}}}}
\makeatother

\DeclareMathOperator*{\argmax}{argmax}

\def\BibTeX{{\rm B\kern-.05em{\sc i\kern-.025em b}\kern-.08em
    T\kern-.1667em\lower.7ex\hbox{E}\kern-.125emX}}
\begin{document}

\title{Efficient Gradient-Based Optimization for Joint Layout Design and Control of Wind Turbines\\
\thanks{ Supported by the Energy Earthshots Initiative as part of the DOE Office of Biological and Environmental Research. PNNL is operated by DOE by the Battelle Memorial Institute under Contract DE-A06-76RLO 1830.}
}

\author{\IEEEauthorblockN{James Kotary}
\IEEEauthorblockA{\textit{Pacific Northwest National Laboratory}\\
Richland, WA \\
james.kotary@pnnl.gov}
\and
\IEEEauthorblockN{Natalie Isenberg}
\IEEEauthorblockA{\textit{Pacific Northwest National Laboratory}\\
Richland, WA \\
natalie.isenberg@pnnl.gov}
\and
\IEEEauthorblockN{Draguna Vrabie}
\IEEEauthorblockA{\textit{Pacific Northwest National Laboratory}\\
Richland, WA \\
draguna.vrabie@pnnl.gov}
}

\maketitle
\begin{abstract}
A central challenge in the design of energy-efficient wind farms is the presence of wake effects between turbines. When a wind turbine harvests energy from free wind, it produces a turbulent region with reduced energy for downstream turbines. Strategies for increasing the efficiency of wind farms by mitigating wake effects have been the subject of much computational research. To this end, both the positioning of the turbines and their cooperative control must be taken into account. Since they are mutually dependent, increasing attention has been paid to the joint consideration of layout design and control. However, joint modeling approaches lead to large-scale and difficult optimization problems, which lack efficient solution strategies. This paper describes an efficient optimization framework for the joint design and steady-state control of wind turbines under a continuous wake effect model. Through a multi-level reformulation of the joint optimization problem, it shows how problem sub-structure can be efficiently exploited to yield order-of-magnitude reduction in convergence times over naive approaches and prior proposals. 
\end{abstract}
\section{Introduction}
Wind is considered one of the primary sources of renewable energy, along with solar and hydropower. In the year $2024$, over $8.1$ percent of global electricity was supplied by wind, and onshore wind in particular has become one of the least expensive energy sources available \cite{abdelilah2020renewables}. On the other hand, the energy productivity of a wind farm is highly dependent on its efficient design and control. A primary source of energy loss in wind farms comes due to aerodynamic interactions between turbines known as wake effects - when free wind energy is collected by a turbine, a turbulent region forms, leaving less energy available to further downstream turbines in its wake. 

A significant body of research is dedicated to optimization models and methods aimed at minimizing the impact of wake effects on power production, via either its control strategies \cite{boersma2017tutorial} or its layout \cite{pranupa2023review}. For example, rotational control around the yaw axis allows for the redirection of wakes away from downstream turbines. Similarly, an optimal layout can position turbines to avoid the worst of their upstream counterparts' wake effects. Layout optimization has widely been considered in isolation from control considerations, but such an approach inherently makes simplifying assumptions about the subsequent control policy. Recent work \cite{chen2022joint} shows that increased power efficiency can be obtained by optimizing the wind farm's layout \emph{jointly} with its steady-state control. However the requisite joint layout and steady-state control problems pose significant computational challenges - they inherit the nonconvexity of the underlying layout problem, while variable discretization over control scenarios leads to impractical problem sizes for conventional solvers.

To address that challenge, this paper presents a gradient-based optimization methodology for efficiently solving the joint wind farm layout and control problem. A specialized algorithm is proposed, based on repeated solution of efficiently solvable control subproblems, whose relation to a master layout problem is efficiently captured by well-defined gradients. Its efficiency advantage is demonstrated on an array of challenging test cases, which show orders-of-magnitude faster convergence over naive approaches on large problems. An open source optimization library comes with this paper, available at: \url{https://github.com/anon-research-coding/WFLCopt}.

\section{Wind Farm Power Generation Model}
This section describes the power production of a wind farm as a function of layout and control variables. We adopt a steady-state parametric wake effect model based on conservation of momentum which is widely used in layout optimization \cite{park2015layout,guirguis2016toward,chen2022joint}. Let there be $N_T$ turbines each of diameter $D$, to be placed within a feasible farm region  $\mathcal{F} \subset \mathbb{R}^2$. The \emph{layout} of the wind farm is defined by positional variables $\bm{l} \in \mathbf{R}^{N_T \times 2}$, within their joint feasible set $\mathcal{L} \coloneqq \mathcal{F}^{N_T}$. Turbines are primarily \emph{controlled} via their yaw angles $\bm{\lambda} \in \mathbb{R}^{N_T}$  within feasible limits $\bm{\lambda} \in \Lambda$, which take the form of upper and lower bounds $\Lambda = \{\bm{\lambda} \textsc{\;\;s.t.\;\;} \lambda_{min} \leq \bm{\lambda} \leq \lambda_{max} \}$ on the rotational angle of each turbine. The layout of a wind farm is considered to be fixed once decided, while the controls allow for adaptation in response to changing wind conditions. We follow the typical modeling approach where wind conditions and associated control states are modeled in the steady state, meaning that dynamic control between states is not explicitly  considered.  

The \emph{power} generated by the wind farm under \emph{freestream wind direction} $\theta^{\omega}$, given positions $\bm{l}$ and yaws $\bm{\lambda}$, is denoted $P( \bm{l},  \bm{\lambda},  \theta^{\omega})$. Its value depends on the effective wind speed at each turbine location, which must account for wake effects from all other turbines. For notational convenience we specify two vectors of independent coordinates $\bm{x}, \bm{y} \in \mathbf{R}^{N_T}$ so that $\bm{l} = \left[\bm{x}; \bm{y}\right] \subset \mathbf{R}^{N_T \times 2}$ is made up of their concatenation. We define additional coordinate systems whose canonical directions are aligned with, and perpendicular to a fixed wind angle $\theta^{\omega}$. For any two turbines indexed $i$ and $j$, we denote their relative angles in Cartesian space as $\theta_{ij}$ and define their inter-distances along the \emph{downstream} and \emph{radial} directions respectively: 

\begin{minipage}{0.23\textwidth}
\vspace{-10pt}
\begin{equation*}
\label{eq:downstream_coord}
         d_{ij} = \| l_i - l_j  \|_2 \cos(\theta_{ij} - \theta^{\omega}),
\end{equation*}
\end{minipage}
\hfill
\begin{minipage}{0.23\textwidth}
\vspace{-10pt}
\begin{equation*}
    \label{eq:radial_coord}
         r_{ij} = \| l_i - l_j  \|_2 \sin(\theta_{ij} - \theta^{\omega}),
\end{equation*}
\end{minipage}

\noindent
where dependence of $d_{ij}$ and $r_{ij}$ on the wind angle index $\omega$ is omitted from their notation, being apparent from context throughout. Corresponding matrix variables $\bm{d}, \bm{r} \in \mathbb{R}^{N_T \times N_T}$ are made up of those scalar values. They represent relative coordinates, induced by the wind direction and centered at each turbine, and serve simplify the modeling of wake effects. 

\subsection{Wake Effect Model}
\label{sec:wake_effect_model}

Wake effects are measured in terms of their resulting windspeed deficits, relative to the \emph{freestream wind speed} $U$.  At each point in $2D$ space, wake effects from upstream turbines propagate primarily along the downstream direction $\bm{d}$, with their deficit magnitudes decaying along the radial direction $\bm{r}$ in a Gaussian profile. Let $WD(\bm{d},\bm{r})$ represent the pairwise \emph{windspeed deficit factors} due to wake effects between turbines: 
\begin{equation}
\label{eq:windspeed_deficit}
WD(\bm{d},\bm{r})  =  \left(1 - \sqrt{1 - \frac{\sigma_{\bm{r}0}}{\sigma_{\bm{r}}}\bm{C_T}}\;   \right)  \exp \left( \frac{(\bm{r} - \bm{\delta}_f)^2}{2\bm{\sigma}_{\bm{r}}^2}  \right),   
\end{equation}
where $\bm{C_T}$ are thrust coefficients. These define Gaussian wake profiles which expand linearly   in the downstream direction at rate $k_{\bm{r}}$. That is, $\bm{\sigma}_{\bm{r}} = \bm{\sigma}_{\bm{r}0} + k_{\bm{r}}\bm{d}$ where $\bm{\sigma}_{\bm{r}0} = \frac{D}{2\sqrt{2}}\cos(\bm{\lambda})$.

Importantly, wind turbines can only induce wake effects on downstream positions. Correspondingly, formula \eqref{eq:windspeed_deficit} holds only for elements where $d_{ij}>0$; otherwise, the deficit factor is typically defined to be zero \cite{park2015layout,chen2022joint}. This leads to potential discontinuities in the function \eqref{eq:windspeed_deficit} which are physically unrealistic and numerically challenging. These can be closely approximated everywhere by a smooth function $\widetilde{WD}(\bm{d},\bm{r}) \coloneqq s(\bm{d}) \cdot WD(\bm{d},\bm{r})$, where $s(\bm{d}) = \frac{1}{1+e^{-\tau \bm{d}}}$ is a sigmoid function, noting that $WD=\widetilde{WD}$ as $\tau \to \infty$. In addition to its linear expansion, the wake profile also experiences deflection of its centerline by the value $\bm{\delta}_f$, also a function of $\bm{d}$:
\begin{subequations}
    \label{eq:wake_deflection}
    \begin{align*}
        \bm{\delta}_f = &\; a_d \cdot D + b_d \cdot \bm{d}  \\ & + \frac{\bm{\phi}}{5.2} E_0 \cdot \sqrt{\frac{\sigma_{\bm{r}0}}{k_{\bm{r}}C_T}} \cdot \ln \left[ \frac{(1.6+\sqrt{C_T})(1.6\sqrt{\nicefrac{\sigma_{\bm{r}}}{\sigma_{\bm{r}0}}}-\sqrt{C_T})}{(1.6-\sqrt{C_T})(1.6\sqrt{\nicefrac{\sigma_{\bm{r}}}{\sigma_{\bm{r}0}}}+\sqrt{C_T})}  \right] 
    \end{align*}
\end{subequations}
where we define intermediate variables: $\bm{E}_0 = \bm{C}_0^2 - 3e^{\nicefrac{1}{12}}C_0 + 3e^{\nicefrac{1}{3}}$ and $\bm{C}_0 = 1 - \sqrt{1 - \bm{C}_T}$, plus $\bm{\phi} = \frac{0.3\bm{\lambda}}{\cos \bm{\lambda}} (1 - \sqrt{1 - \bm{C}_T \cos \bm{\lambda}})$ following \cite{chen2022joint}. Detailed formulations of the deflection model are found in \cite{bastankhah2016experimental}, and $a_d$, $b_d$ and $k_{\bm{r}}$ are physical parameters.

\subsection{Effective Wind Speed and Power Production}
The amount of power that can be harvested by turbine $i$ is determined by the \emph{effective windspeed} $V_i$ at its location, which results from the \emph{combined} wake effects from all turbines in the wind farm. To combine multiple wake effects at a single point, we adopt the traditional approach which determines the total windspeed deficit as the Euclidean norm over those induced by the individual turbines \cite{park2015layout,guirguis2016toward,chen2022joint}:
\begin{equation}
    \label{eq:effective_windspeed}
         V_i  =  U \cdot \left(1 - \sqrt{\sum_{j \neq i} \widetilde{WD}(d_{ji},r_{ji})^2}  \right).   
\end{equation}

The total power generated by the wind farm is the sum over that of each turbine, due to its effective wind speed:
\begin{equation}
    \label{eq:power_fn}
        P(\bm{l}, \bm{\lambda},  \theta^{\omega}) = \sum_{i=1}^N \frac{1}{2} \rho (\frac{\pi}{4}D^2) C_{P_i} \cos (\lambda_i) V_i^3.
\end{equation}
The coefficients $\bm{C_P}$, along with $\bm{C_T}$ in \eqref{eq:windspeed_deficit}, depend on controllable axial induction factors $\alpha_i$ for each turbine. While nothing prevents their optimization along with $\bm{\lambda}$, they are considered less effective than yaw control \cite{van2019effects,fleming2017field} and often set to the constant value $\alpha_i = \nicefrac{1}{3}$, which we adopt.

Readers interested in further details and derivations of the power generation models are referred to \cite{park2015layout,bastankhah2016experimental}. The goal of this work is optimization of power objectives based on \eqref{eq:power_fn}, for which it suffices to recognize that \eqref{eq:power_fn} is a smooth but highly nonconvex function of layout and control variables $(\bm{l}, \bm{\lambda})$.
\section{Joint Optimization Problem}
\label{sec:joint_problem}
The aim of this paper is to develop optimization methods for solving a \emph{joint} layout and control problem, which can be viewed as an optimization problem in two stages. The first-stage decision is permanent and determines the layout $\bm{l} \in \mathcal{L}$, while the second determines control variables $\bm{\lambda} \in \Lambda$. Their roles are distinguished by the fact that $\bm{\lambda}$ can be adjusted in response to changing wind conditions $\theta^{\omega}$. To quantify those changing conditions, we assume that $\theta^{\omega}$ are realized from an underlying continuous random variable $\Theta$. Like prior works \cite{park2015layout,guirguis2016toward,chen2022joint}, we approximate $\Theta$ by deriving from it a discrete random variable in order to form a tractable optimization problem. First, the range $[0^o, 360^o]$ is partitioned into $W$ equal-sized bins. Each $\omega^{th}$ bin is represented by its midpoint, denoted $\theta^{\omega}$ for $1 \leq \omega \leq W$. These wind directions define $W$ distinct \emph{scenarios}, each in which \emph{scenario-wise turbine control variables} $\bm{\lambda}^{\omega}$ may be optimized in response. Wind directions $\theta^{\omega}$ are sampled from $\Theta$, and the frequency of each bin is used to construct the associated scenario probabilities $p^{\omega}$. All wind directions are relative to a fixed reference angle at due East.

The overall \emph{joint layout and control problem} seeks a single assigment of turbine positions $\bm{l}$, plus yaw angles $\bm{\lambda}^{\omega}$ for each of $W$ wind directions $\theta^{\omega}$, which maximize the expected power:
\begin{subequations}
\label{eq:joint_opt}
\vspace{-12pt}
\begin{align}
     \max_{\bm{l}, \bm{\lambda}^{\omega}   \; \forall \omega} & \;\; \sum_{\omega = 1 }^W p^{\omega} P(\bm{l}, \bm{\lambda}^{\omega},  \theta^{\omega}) \label{eq:joint_opt_obj} \\
    \textit{s.t.} \;\;\;& \;\; \bm{l}\in \mathcal{L} \label{eq:joint_constraint_l} \\
     \;\;\;& \;\; \bm{\lambda}^{\omega} \in \Lambda \;\; \;\;\;\;\;\;\;\;\;\;\;\; \forall \omega \;\;\;\; \textit{s.t.}  \;\;  1 \leq \omega \leq W \label{eq:joint_constraint_lamb}  \\
    \;\;\;& \;\;\| \bm{l}_i - \bm{l}_j \|^2_2 \geq 4D \;\;\;\;\forall i,j\;\;  \textit{s.t.} \;\; 1\leq i<j \leq N_T \label{eq:joint_opt_interdistance_constr}.
\end{align}
\end{subequations}
In addition to the wind farm perimeter \eqref{eq:joint_constraint_l} and control bounds \eqref{eq:joint_constraint_lamb},  inter-distance constraints \eqref{eq:joint_opt_interdistance_constr} ensure that wind turbines are spaced at a minimum safe distance. We summarize the main challenges posed by problem \eqref{eq:joint_opt} as follows: \textbf{(1)} Its variable count scales quadratically as $N_T \cdot (W+2)$, quickly leading to problem sizes which challenge conventional solvers and  \textbf{(2)} It is nonconvex in both objective and constraints, with flat regions and sharp inclines prevalent in each term of \eqref{eq:joint_opt_obj}, while nonconvex inter-distance constraints number quadratically in $N_T$. On the other hand, the majority of variables in \eqref{eq:joint_opt} are scenario-wise variables $\bm{\lambda}^{\omega}$, which are mutually uncoupled and subject only to simple bounds. The following proposal leverages the special structure in Problem \eqref{eq:joint_opt} to arrive at a framework for its efficient optimization.
\vspace{-2pt}
\section{Multilevel Reformulation and Solution}
\label{sec:MLR}
\vspace{-2pt}
In this section, we present a computational framework for building efficient optimization methods to solve problem \eqref{eq:joint_opt}. We first observe that the optimal control variables $\bm{\lambda}^{\omega}$ in each scenario are directly determined by the layout $\bm{l}$, and the corresponding power production can be differentiated as a function of $\bm{l}$. The \emph{optimal power function} $P^{\star}$ below represents the available power as a function of position, given that yaws are optimized with respect to the wind direction $\theta^{\omega}$:
\begin{equation}
\label{eq:yaw_optimal_power}
      P^{\star}(\bm{l}, \theta^{\omega}) \coloneqq \max_{ \bm{\lambda} \in \Lambda }  \;\;  P(  \bm{l},  \bm{\lambda},   \theta^{\omega}).
\end{equation}

The joint problem \eqref{eq:joint_opt} is then equivalent to the following:
\begin{subequations}
\label{eq:decomp_opt}
\begin{align}
\label{eq:decomp_opt_objective}
     \max_{\bm{l} \in \mathbb{R}^{2 \cdot N_T} } & \;\; \sum_{\omega =1}^W p^{\omega} P^{\star}(\bm{l}, \theta^{\omega}) \\
    \textit{s.t.} \;\;\;& \;\; \bm{l}\in \mathcal{L} \\
    \;\;\;& \;\;\| \bm{l}_i - \bm{l}_j \|^2_2 \geq 4D \;\;\;\;\forall i,j\;  \textit{s.t.} \;\; 1\leq i<j \leq N_T.
\end{align}
\end{subequations}
Problems \eqref{eq:yaw_optimal_power},\eqref{eq:decomp_opt} can be viewed as the \emph{inner} and \emph{outer problems} in a  \emph{multilevel reformulation} of problem \eqref{eq:joint_opt}.  While equivalent to the original joint problem \eqref{eq:joint_opt}, it poses different challenges. Each of its objective terms $P^{\star}(\bm{l}, \theta^{\omega})$ requires (nonconvex) optimization to evaluate, making the overall objective function  expensive to evaluate. The lack of a closed form for each $P^{\star}(\bm{l}, \theta^{\omega})$ also prevents direct application of gradient-based optimization solvers to \eqref{eq:decomp_opt}. 

On the other hand, this formulation is expressed solely in terms of the $2N_T$ positional variables $\bm{l}$, which number a factor of $W$ less than those of the joint formulation. The present proposal is based on directly solving the smaller-scale problem \eqref{eq:decomp_opt} in place of the original joint problem \eqref{eq:joint_opt}, while at the same time mitigating the computational challenges of its objective function. It  relies on two main concepts: \textbf{(1)} exact differentiation of $P^{\star}$ via envelope theorems and \textbf{(2)} an efficient solution scheme for repeated resolution of the inner problems \eqref{eq:yaw_optimal_power}. A computational pipeline composed of these techniques enables fast, repeated evaluation of the composite objective \eqref{eq:decomp_opt_objective} and its gradients, which can be paired with various appropriate nonconvex solvers  to directly optimize problem \eqref{eq:decomp_opt}. Each are detailed in the following subsections.

\subsection{Differentiation of the Implicit Objective Function}
\label{sec:diff_opt_danskins}
Let $f(\bm{l}) \coloneqq \sum_{\omega = 1}^W p^{\omega} \cdot P^{\star}(\bm{l}, \theta^{\omega})$ denote the objective function in problem \eqref{eq:decomp_opt}. Its gradient is found by combining the gradients of $W$ implicitly defined optimal power functions:
\begin{equation}
\label{eq:total_gradient}
    \nabla_{\bm{l}} f(\bm{l}) =  \sum_{\omega = 1}^W p^{\omega} \cdot  \nabla_{\bm{l}} P^{\star}(\bm{l}, \theta^{\omega}) 
\end{equation}
An efficient formula for each constituent term $\nabla_{\bm{l}} P^{\star}(\bm{l}, \theta^{\omega})$ is provided by the Envelope Theorems for differentiation of a continuous minimizer \cite{milgrom2002envelope}, which reduce its derivative to a standard partial derivative of the underlying power function $P$, evaluated at its corresponding argmin. Define $\bm{\lambda}^{\star} = \argmax_{ \bm{\lambda} \in \Lambda }   P(  \bm{l},  \bm{\lambda},   \theta^{\omega})$ so that $\bm{\lambda}^{\star} = \argmax_{ \bm{\lambda} \in \Lambda }   P(  \bm{l},  \bm{\lambda},   \theta^{\omega})$, then
\begin{equation}
\label{eq:danskins}
    \nabla_{\bm{l}} P^{\star}(\bm{l}, \theta^{\omega}) = \nabla_{\bm{l}}   P( \bm{l}, \bm{\lambda}^{\star}, \theta^{\omega}).  
\end{equation}
That is, to compute the gradient $\nabla_{\bm{l}} P^{\star}(\bm{l}, \theta^{\omega})$, one first solves the inner problem \eqref{eq:yaw_optimal_power} for $P^{\star}(\bm{l}, \theta^{\omega})$. The associated argmin $\bm{\lambda}^{\star}$ comes as a byproduct, at \emph{no additional cost}. The remaining standard gradient $ \nabla_{\bm{l}} P( \bm{l}, \bm{\lambda}^{\star}, \theta^{\omega})$ in \eqref{eq:danskins} is easily computed, e.g., by automatic differentiation (see details in Section \ref{sec:Experiments}).

Note that $\nabla_{\bm{l}} P^{\star}(\bm{l}, \theta^{\omega})$ could be approximated by finite differences \cite{froese2022optimal}, which is the default behavior of some black-box optimization solvers. However, this requires solving problem \eqref{eq:yaw_optimal_power} once for each of its $N_T$ components. Thus, approximation of $\nabla_{\bm{l}} f(\bm{l})$ requires $N_T^2$ solver calls. By comparison, this result provides exact gradients with no additional solver calls. 

\smallskip
\emph{Theoretical Remarks:} We remark that the validity of \eqref{eq:danskins} requires $P$ to be continuous and have a unique partial minimizer $\lambda^{\star}$ \cite{milgrom2002envelope}. Continuity of $P$ follows from its construction as a composition of continuous functions \cite{park2015layout}, but nonunique $\lambda^{\star}$ may arise since total power is unaffected by the yaw angles of the most-downstream turbines, since they induce no wake effects. This issue is resolved by adding a tie-breaking term to the power objective, which has negligible effects on the solution but ensures uniqueness. While nonconvexity of $P$ prevents guarantees that $\nabla_{\bm{l}} P^{\star}(\bm{l}, \theta^{\omega})$ exists everywhere, it is typical in nonconvex settings to prove that it exists almost-everywhere, and that formula \eqref{eq:danskins} is valid whenever it exists \cite{milgrom2002envelope}. A detailed analysis on the existence of the gradient \eqref{eq:danskins} is reserved for an extension of this work. 

\subsection{Fast Solution and Updating Scheme for the Inner Problems}
\label{sec:lower_level_sol}
A generic optimization algorithm repeatedly updates an estimate of the decision variable $\bm{l}$, as a deterministic function $U$ of the current estimate, plus its corresponding objective value $f(\bm{l})$ and gradient $\nabla_{\bm{l}} f(\bm{l})$. For this we define notation:
\begin{equation}
\label{eq:solution_update}
    \bm{l} \coloneqq U\left(\; \bm{l}, f(\bm{l}), \nabla_{\bm{l}} f(\bm{l})\; \right).
\end{equation}
Recall that the multilevel reformulation \eqref{eq:decomp_opt} gains a drastic reduction in problem size from the joint problem \eqref{eq:joint_opt}, but at the cost of complicating the evaluation of the arguments $f(\bm{l})$ and $\nabla_{\bm{l}} f(\bm{l})$ to the update routine \eqref{eq:solution_update}. To alleviate that tradeoff, this section thus presents a computational pipeline that provides fast, repeated evaluations of the optimal power functions $P^{\star}(\bm{l}, \theta^{\omega})$ and thus also the total objective $f(\bm{l})$ in problem \eqref{eq:decomp_opt}.  In principle, the framework is solver-agnostic and can be paired with any optimization method represented by \eqref{eq:solution_update} which accomodates nonconvex objectives and constraints.

A major advantage to the reformulation \eqref{eq:decomp_opt} is that its inner subproblems \eqref{eq:yaw_optimal_power}  require optimization subject only to simple bounds, a.k.a. box constraints $\Lambda = \{ \bm{\lambda}: \lambda_{\min} \leq \bm{\lambda} \leq \lambda_{\max} \}$. This overall structure, which has nonconvex box-constrained subproblems within an outer problem subject to arbitrary constraints, is reminiscent of the dual problem formulations employed by state-of-the-art augmented Lagrangian solvers such as LANCELOT \cite{conn2013lancelot}. Their success relies on their efficient handling of the inner problem using a limited-memory box-constrained Broyden–Fletcher–Goldfarb–Shanno (L-BFGS-B) algorithm \cite{liu1989limited}. This algorithm is an ideal choice for two main reasons: \textbf{(1)} It is proven to excel in optimizing nonconvex functions over box constraints and \textbf{(2)} it makes highly efficient use of hot-starts when repeated solutions are desired \cite{conn2013lancelot}.

\begin{figure}[] \includegraphics[width=1.00\columnwidth]{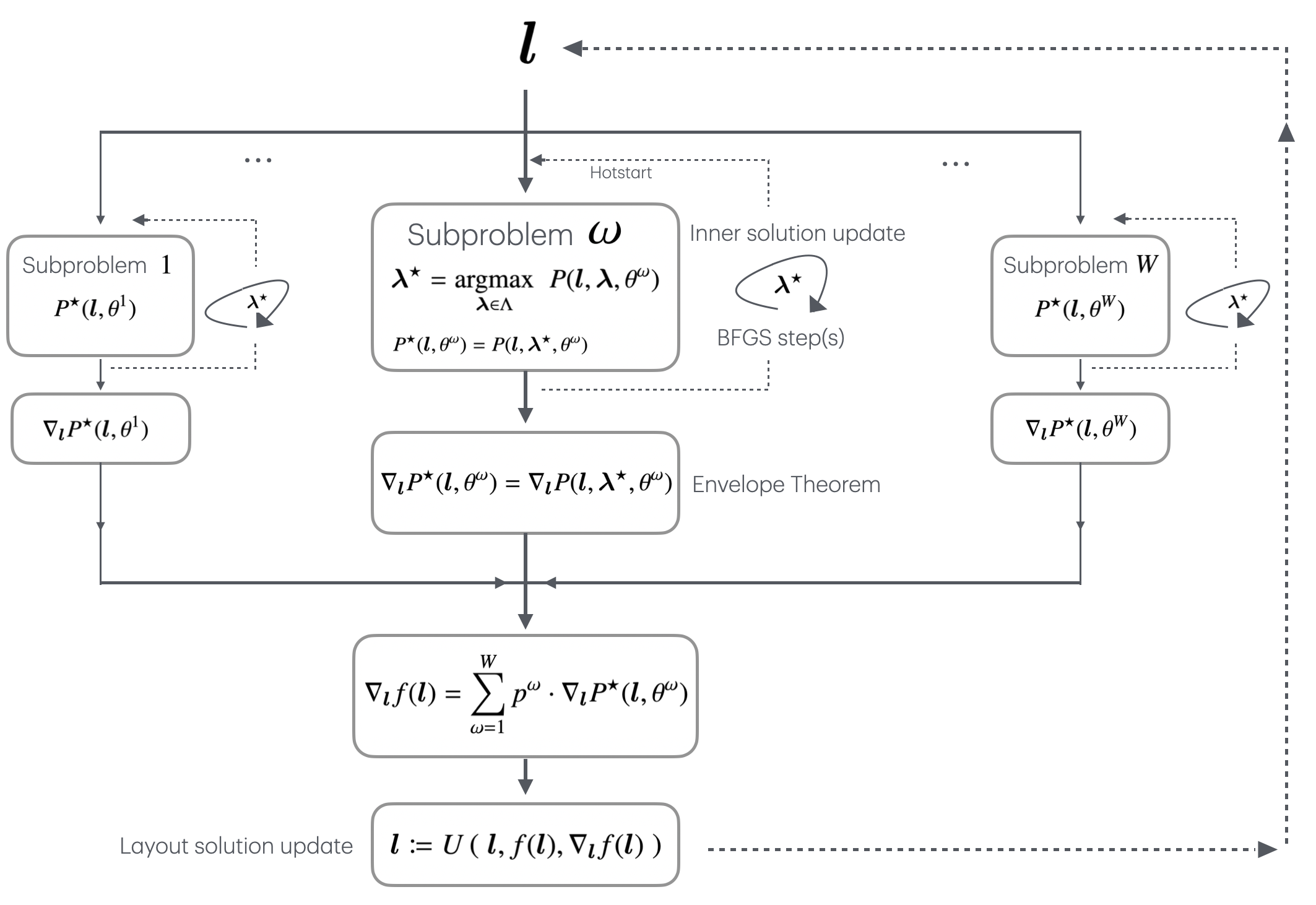}
    \caption{Illustrating the combined computational framework.} 
    \label{fig:illustration}
\vspace{-16pt}
\end{figure}

Our framework applies a similar strategy, but to \emph{each of the constituent subproblems} \eqref{eq:yaw_optimal_power} within the multilevel reformulation \eqref{eq:decomp_opt}. An overall schematic is illustrated in Figure \ref{fig:illustration}. Each iteration of a generic optimization step \eqref{eq:solution_update} revises location variables $\bm{l}$ in the outer problem \eqref{eq:decomp_opt}. To prepare the next iteration, re-evaluation of $f(\bm{l})$ requires re-solving the $W$ inner problems \eqref{eq:yaw_optimal_power} to get $P^{\star}(\bm{l}, \theta^{\omega})$ for $1 \leq  \omega \leq W$. Each re-solve is accelerated by hotstarting the BFGS algorithm from a cached copy its optimal solution $\left( \bm{\lambda}^{\star} \right)^{\omega}$ from the previous outer iteration. Since the outer iterations \eqref{eq:solution_update} result in small changes to $\bm{l}$, this reduces the cost of successive calls to $f(\bm{l})$  by over 99$\%$ after the first few outer iterations in our experiments. Following Section \ref{sec:diff_opt_danskins}, $\nabla_{\bm{l}} f(\bm{l})$ is then efficiently obtained by computing the standard derivatives \eqref{eq:danskins} with respect to each subproblems and combining them per equation \eqref{eq:total_gradient}. The next section shows empirically how it enables order-of-magnitude reductions in solving time when paired with off-the-shelf optimization solvers to implement the iterations \eqref{eq:solution_update}.
\section{Implementation and Experiments}
\label{sec:Experiments}
{\renewcommand{\arraystretch}{1.0}
\begin{table*}[ht]
\small
\centering
\caption{Algorithm Runtimes and Objective Values Across Problem Sizes}
\label{tab:Results}
\resizebox{\textwidth}{!}{%
\begin{tabular}{ll*{12}{c}}
\toprule
\multirow{3}{*}{\textbf{Form}} & 
\multirow{3}{*}{\textbf{Method}} & 
\multicolumn{12}{c}{\textbf{Number of Turbines $N_T$}} \\
\cmidrule(lr{0.3em}){3-14}
& & 
\multicolumn{6}{c}{\textbf{Runtime (s)}} & 
\multicolumn{6}{c}{\textbf{Objective Value (GWh)}} \\
\cmidrule(lr{0.3em}){3-8} \cmidrule(lr{0.3em}){9-14}
& & 
\textbf{25} & \textbf{36} & \textbf{49} & \textbf{64} & \textbf{81} & \textbf{100} & 
\textbf{25} & \textbf{36} & \textbf{49} & \textbf{64} & \textbf{81} & \textbf{100} \\
\midrule
\multirow{2}{*}{JOINT} 
  & IP  & 3719 & 10147 & 29288 & 57565 & \ding{55} & \ding{55} & 475.1 & 672.9 & 899.9 & 1157.5 & \ding{55} & \ding{55}  \\
  & SQP & 1225 & 4315 & 15363 & 52469 & 109790 & \ding{55} & \textbf{476.9} & 674.1 & 903.3 & 1163.0 & 1457.2 & \ding{55}  \\
\midrule
\multirow{2}{*}{\textbf{MLR-A}} 
  & IP  & 3404 & 4780 & 7414 & 14795 & 27837 & 55954 & 475.6 & 673.5 & 902.2 & 1162.2 & 1454.2 & 1773.6  \\
  & \textbf{SQP} & \textbf{364} & \textbf{569} & \textbf{933} & \textbf{1873} & \textbf{3353} & \textbf{6976} & \textbf{476.9} & \textbf{674.8} & \textbf{904.5} & \textbf{1165.6} & \textbf{1457.9} & \textbf{1780.3}  \\
\bottomrule
\end{tabular}%
}
\end{table*}
}
In this section, we evaluate the ability of our proposed computational framework to improve upon the efficiency of conventional optimization solvers applied to the joint layout and control problem \eqref{eq:joint_opt}. Section \ref{sec:methods} discusses implementation details regarding each method evaluated. An array of increasingly difficult test problems is described in Section \ref{sec:problem_specs}, and a systematic comparison of results is presented in \ref{sec:results}. 

\begin{table}[h!]
\centering
\caption{Physical Parameters}
\begin{tabular}{c c c c c c c c c c}
\hline
$D$ & $\rho$ & $U$ & $a_d$ & $b_d$ & $k_{\bm{r}}$ & $\lambda_{\min}$ & $\lambda_{\max}$ \\
\hline
126m & $1.23\nicefrac{kg m}{s}$ & $8\nicefrac{m}{s}$ & -0.035 & -0.01 &0.03& $-30^{o}$ & $30^{o}$ \\
\hline
\end{tabular}
\vspace{-8pt}
\end{table}
\subsection{Methods and Implementation}
\label{sec:methods}
The computational framework described in Section \ref{sec:MLR} is independent of any particular optimization algorithm used to solve the outer problem \eqref{eq:decomp_opt}, which so far has been generically represented by equation \eqref{eq:solution_update}. When considering the choice of algorithm, it must be capable of handling general nonlinear programs with nonconvex constraints and objective. It should also be equipped with robust stepsize control, including advanced heuristics such as trust region and line search methods. We consider two main classes of algorithms which share these properties: \textbf{(1)} the \emph{interior point methods} (IP), where $U$ in \eqref{eq:solution_update} solves a log-barrier approximation to the true problem at each step and \textbf{(2)} \emph{sequential quadratic programming (SQP)}, where $U$ solves a quadratic approximation to the problem. We refer to these as the \emph{base methods}. We apply them to each \emph{problem form}: \textbf{(1)} the original joint form \eqref{eq:joint_opt} (called JOINT), and \textbf{(2)} the multilevel reformulation \eqref{eq:decomp_opt} augmented by the techniques proposed in Section \eqref{sec:MLR} (called MLR-A). Each problem form / base method combination is evaluated, totaling four. The JOINT variants  \emph{serve as baselines} in this study, having been previously proposed for isolated layout optimization \cite{park2015layout,guirguis2016toward} and tested on the joint layout and control problem \cite{chen2022joint}.

\emph{Implementation:} All optimization methods tested in this section are implemented using the SciPy suite of  solvers. Its  interior point method is a trust-region IP variant known as NITRO \cite{byrd1999interior}. Its main SQP method is a variant called SLSQP  \cite{kraft1988software}. Finally, our MLR-A methods are built using SciPy's implementation of L-BFGS-B, based on the proposal of \cite{liu1989limited}. All solvers are given an objective tolerance of $1e-5$, which creates a fair standard for comparing convergence times between equivalent objectives \eqref{eq:joint_opt_obj} and \eqref{eq:decomp_opt_objective}.  All other settings use SciPy's standard values. Derivatives are calculated by automatic differentiation in PyTorch, including \eqref{eq:danskins}.

\emph{On Further Baselines:} We remark that prior work \cite{chen2022joint} proposes a decomposition scheme for solving \eqref{eq:joint_opt}, but \emph{does not report speedups} over naive SQP on the JOINT formulation. The included library contains our Python implementation of \cite{chen2022joint}. Due to this lack of runtime advantage, plus inconsistent convergence in our studies, a systematic comparison with \cite{chen2022joint} is reserved for an extension of this work.

\subsection{Problem Specifications}
\label{sec:problem_specs}
To evaluate the scalability of each method, we solve a sequence of joint layout and control problems in which $N_T \in \{ 25, 36, 49, 64, 81, 100 \}$ increases by perfect squares. In each case, $W=36$ is held constant, and $\theta^{\omega}$ are sampled for $ 1 \leq \omega \leq W$ from a Dirichlet distribution with $\alpha=1$. Feasible turbine positions are in $ \cF =  \left[0,x_{max}\right] \times \left[0,y_{max}\right]$. We set $x_{max} = y_{max} \coloneqq {1.3} \sqrt{N_T} \cdot 4D$ to standardize the turbine density across cases. Problems \eqref{eq:joint_opt} and \eqref{eq:decomp_opt} are prone to local optima, which have a randomizing effect on both the convergence time and final objective value. Thus, \emph{all Runtime and Objective values are averaged over $30$ independent, randomly initialized solver runs}.

\section{Results and Conclusions}
\label{sec:results}
\begin{figure}[] \includegraphics[width=1.0\columnwidth]{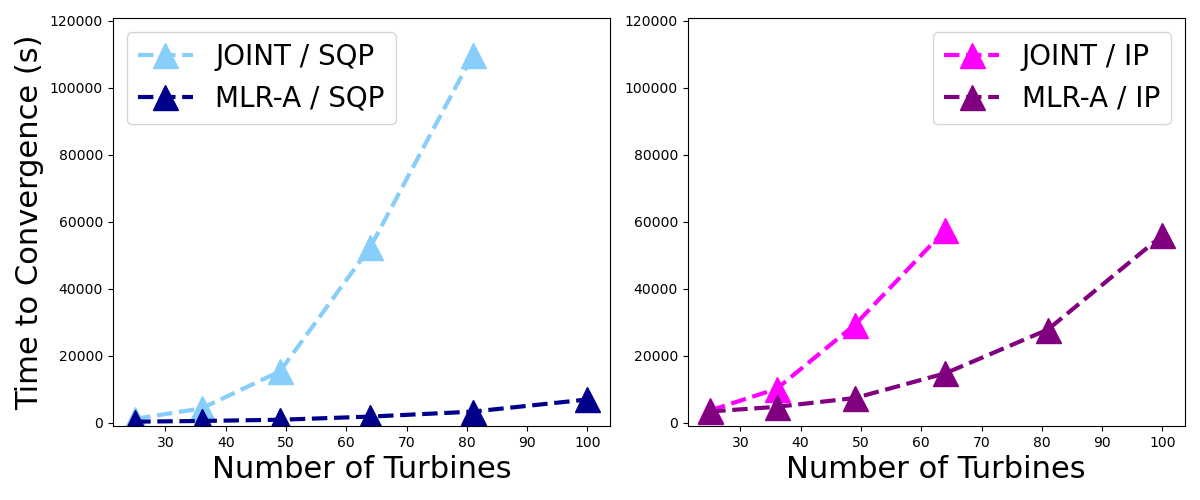}
    \caption{Runtime vs Problem Size: SQP (left) and IP (right) variants, applied directly to problem \eqref{eq:joint_opt} vs MLR framework.} 
    \label{fig:timing}
\vspace{-20pt}
\end{figure}

Table \ref{tab:Results} reports the Runtime of each solution in seconds, along with the Objective Value in GWh, after multiplying the nominal power objective  \eqref{eq:joint_opt_obj},\eqref{eq:decomp_opt_objective} by $T=8760 \nicefrac{hr}{yr}$ to yield annual energy. The best results across each combination of problem form and base method are marked in bold. Runtimes are further illustrated in Figure \ref{fig:timing} as a function of problem size. A \emph{time limit} of $36$hr is imposed on each run. Runs that exceed that limit are marked by \ding{55} and omitted from the Figure. 

Figure \ref{fig:timing} shows how both the SQP (at left) and IP (at right) variants of MLR-A gain an increasing runtime advantage over their baseline JOINT counterparts. Both IP and SQP exceed the timeout threshold on the JOINT problem, past $N_T=64$ and $N_T=81$ respectively. Meanwhile the MLR-A approach converges on all problem sizes when paired with both the IP and SQP base methods. Regarding objective value, the MLR-A approach also consistently reports a marginal advantage, giving best results when paired with SQP. We speculate that the reduction of variable space leads to fewer local optima, but this effect is secondary to the runtime advantage.

We conclude that MLR-A should be combined with SQP for maximum efficiency. The comparison with IP methods serves as an ablation which supports this finding, which is consistent with previous findings which favor SQP for layout optimization. \emph{Our MLR-A variant (dark blue) converges up to $33$ times faster than the standard SQP approach (light blue) while also gaining in the objective value}.

\smallskip \noindent
\emph{Concluding Remarks:} This study shows the potential for specialized optimization strategies to enhance automated design of wind farm layouts. The resulting runtime reductions may encourage further interest in gradient-based approaches, and open the door to further opportunities for integration with alternatives such as metaheuristics. Generalization to floating offshore and other settings are directions for future work.

\bibliography{lib}
\bibliographystyle{IEEEtran}

\newpage

\end{document}